\documentclass[a4paper,12pt]{article}
\usepackage{amssymb}
\usepackage{amsthm}
\newtheorem{Lemma}{Lemma}
\newtheorem{Proposition}{Proposition}
\newtheorem{Theorem}{Theorem}
\newtheorem{Remark}{Remark}
\newtheorem{Corollary}{Corollary}
\newtheorem{Definition}{Definition}

\newcommand{\Pic}{\mathrm{Pic}}

\newcommand{\codim}{\mathrm{codim}}
\newcommand{\Hilb}{\mathrm{Hilb}}
\newcommand{\hilb}{\mathrm{hilb}}
\newcommand{\mod}{\mathrm{mod.}}

\newcommand{\Q}{\mathbb{Q}}
\newcommand{\C}{\mathbb{C}}
\renewcommand{\P}{\mathbb{P}}
\newcommand{\mg}{\mathcal{M}_g}
\newcommand{\mgo}{\mathcal{M}_g^0}
\newcommand{\mgbar}{\overline{\mathcal{M}}_g}
\newcommand{\pdg}{P_{d,g}}
\newcommand{\pdgbar}{\overline{P}_{d,g}}
\newcommand{\pdC}{\overline{P}_{d,C}}
\newcommand{\pdX}{\overline{P}_{d,X}}
\newcommand{\pdgfun}{\overline{\mathcal{P}}_{d,g}}
\newcommand{\sgbar}{\overline{S}_g}
\newcommand{\sgfun}{\overline{\mathcal{S}}_g}

\title{On the geometry of the compactification of the universal Picard variety}
\author{Claudio Fontanari}
\date{}
\begin{document}
\maketitle
\begin{small}
\begin{center}
\textbf{Abstract}
\end{center}
Here we focus on the geometry of $\pdgbar$, the compactification of the 
universal Picard variety constructed by L.~Caporaso.
In particular, we show that the moduli space of spin curves constructed 
by M.~Cornalba naturally injects into 
$\pdgbar$ and we give generators and relations of the rational divisor class 
group of $\pdgbar$, extending previous work by A.~Kouvidakis.
\end{small}

\section{Introduction} 
The universal Picard variety $\pdg$ is the coarse moduli space for line 
bundles of degree $d$ on smooth algebraic curves of genus $g$. 
Even though one is mainly interested in the behaviour of line bundles 
on smooth curves, nevertheless it is often useful to control their 
degenerations on singular curves. Perhaps the most celebrated example of 
proof by degeneration is provided by the Brill-Noether-Petri theorem 
(see \cite{ACGH:85} and the references therein). Another very recent 
achievement of degeneration techniques is the proof given by 
L.~Caporaso and E.~Sernesi (see \cite{CapSer:00} and \cite{CapSer:02}) 
that a general curve of genus $g \ge 3$ can be recovered from its odd   
theta-characteristics. In particular, in order to control degenerations 
of curves with prescribed theta-characteristics, a key r\^ ole is played 
in \cite{CapSer:02} by the moduli space of spin curves $\sgbar$ 
constructed by M.~Cornalba in \cite{Cornalba:89}. 
This perspective suggests the deep mathematical interest (both in itself 
and as a tool) of a geometrically meaningful compactification of 
the moduli spaces parameterizing pairs of curves and line bundles. 
Let $\pdgbar$ denote the compactification of $\pdg$ constructed by 
L.~Caporaso in \cite{Caporaso:94} via geometric invariant theory. 
The boundary points of $\pdgbar$ correspond to certain line
bundles on  Deligne-Mumford semistable curves, while 
all previously known compactifications of the
generalized Jacobian of an integral nodal curve used torsion free 
sheaves of rank one. From this point of view, a strict analogy 
emerges between $\pdgbar$ and $\sgbar$: even though the techniques 
used in the two constructions are completely different, in both cases 
the resulting compactification is given in terms of line bundles on the
same kind of singular curves. We will see that this analogy has a precise 
explanation: namely, in section~\ref{spin} we introduce a subscheme 
of $\pdgbar$ which compactifies the locus in $\pdg$ corresponding to curves 
with theta characteristics and we investigate how it is related to $\sgbar$. 
Indeed, after having established the existence of a natural morphism 
from $\sgbar$ to $\pdgbar$ (see Theorem~\ref{naturalmorphism}), we prove 
that it is an injection (see Theorem~\ref{injection}) and we give a precise 
combinatorial description of its image (see Theorem~\ref{image}). 
In the future we hope to address similar questions also for moduli spaces of
higher spin curves, which were introduced by T.~J.~Jarvis in a rather 
different style (see \cite{Jarvis:98}, \cite{Jarvis:00}, and \cite{Jarvis:01}).
In section~\ref{divisors}, instead, we obtain a complete description 
of the rational divisor class group of $\pdgbar$ (see Theorem~\ref{A_4g-4}); 
its rational Picard group is determined whenever $\pdgbar$ turns out 
to be a geometric quotient (see Corollary~\ref{Pic}). The strategy of proof is 
straightforward: first of all, we deduce from the basic properties of $\pdgbar$
a rough description of its boundary (see Proposition~\ref{unique}); next, we 
recall a theorem proved by A.~Kouvidakis in \cite{Kouvidakis:91} on the 
Picard group of $\pdg$ (see Theorem~\ref{Kouvidakis}). Hence the result on 
generation follows and in order to exclude nontrivial relations
we simply lift to $\pdgbar$ the families of curves constructed by 
E.~Arbarello and M.~Cornalba in \cite{ArbCor:87}.

We work over the field $\C$ of complex numbers.

It is a great pleasure to thank Enrico Arbarello, who in many fruitful 
conversations shared with me his deep insight and invaluable experience, 
and Lucia Caporaso, who closely followed my work with her contagious 
enthusiasm and patient support.

This research was partially supported by MIUR (Italy).

\section{Notation and preliminaries}\label{notation}
Let $X$ be a Deligne-Mumford semistable curve and let $E$ be 
a complete, irreducible subcurve of $X$. One says that $E$ is 
\emph{exceptional} if it is smooth, rational, and meets the other 
components in exactly two points. Moreover, one says that $X$ is 
\emph{quasistable} if any two distinct exceptional components of $C$ 
are disjoint. In the sequel, $\tilde{X}$ will denote the subcurve 
$\overline{X \setminus \cup E_i}$ obtained from $X$ by removing all  
exceptional components. 

A \emph{spin curve} of genus $g$ (see \cite{Cornalba:89}, \cite{Cornalba:91})
is the datum of a quasistable genus $g$ curve $X$ with an invertible 
sheaf $\zeta_X$ of degree $g-1$ on $X$ and a homomorphism of invertible 
sheaves
$$
\alpha_X: \zeta_X^{\otimes 2} \longrightarrow \omega_X
$$
such that 

(i) $\zeta_X$ has degree $1$ on every exceptional component of $X$;

(ii) $\alpha_X$ is not zero at a general point of every non-exceptional 
component of $X$.

\noindent From the definition it follows that $\alpha_X$ vanishes identically 
on all exceptional components of $X$ and induces an isomorphism 
$$
\tilde{\alpha}_X: \zeta_X^{\otimes 2} \vert_{\tilde{X}} \longrightarrow 
\omega_{\tilde{X}}.
$$ 
In particular, when $X$ is smooth, $\zeta_X$ is just a theta-characteristic 
on $X$. 

\noindent By definition (see \cite{Cornalba:89}, \S~2), two spin curves 
$(X, \zeta_X, \alpha_X)$ and $(X', \zeta_{X'}, \alpha_{X'})$ are 
\emph{isomorphic} if there are isomorphisms $\sigma: X \to X'$ and 
$\tau: \sigma^*(\zeta_X') \to \zeta_X$
such that $\tau$ is compatible with the natural 
isomorphism between $\sigma^*(\omega_{X'})$ and $\omega_X$.
However, we point out the following fact.

\begin{Lemma}~\label{isomorphic}
Let $(X, \zeta_X, \alpha_X)$ and $(X', \zeta_{X'}, \alpha_{X'})$ 
be two spin curves and assume that there are isomorphisms
$\sigma: X \to X'$ and $\tau: \sigma^*(\zeta_{X'}) \to \zeta_X$.
Then $(X, \zeta_X, \alpha_X)$ and $(X', \zeta_{X'}, \alpha_{X'})$
are isomorphic as spin curves.
\end{Lemma}

\proof Let $\beta_X: \zeta_X^{\otimes 2} \to \omega_X$ be defined by the 
following commutative diagram:
$$
\hspace{1.8cm} \zeta_X^{\otimes 2} \hspace{1.3cm} 
\stackrel{\beta_X}{\longrightarrow} 
\hspace{0.7cm} \omega_X
$$
$$
\hspace*{2cm} \uparrow \hspace{0.1cm} \tau^{\otimes 2}
\hspace*{2.5cm} \uparrow \hspace{0.1cm} \cong
$$
$$
\hspace{0.2cm}
(\sigma^*\zeta_{X'})^{\otimes 2} = \sigma^*(\zeta_{X'}^{\otimes 2})
\hspace{0.2cm}
\stackrel{\sigma^* \alpha_{X'}}{\longrightarrow} 
\hspace{0.2cm} \sigma^* \omega_{X'}.
$$
Then $(X, \zeta_X, \beta_X)$ is a spin curve, which is isomorphic to
$(X', \zeta_{X'}, \alpha_{X'})$ by definition and to 
$(X, \zeta_X, \alpha_X)$ by \cite{Cornalba:89}, Lemma~(2.1), so the thesis 
follows.

\qed 

\noindent A \emph{family of spin curves} is a flat family of quasistable 
curves $f: \mathcal{X} \to S$ with an invertible sheaf $\zeta_f$ on 
$\mathcal{X}$ and a homomorphism
$$
\alpha_f: \zeta_f^{\otimes 2} \longrightarrow \omega_f
$$
such that the restriction of these data to any fiber of $f$ gives rise 
to a spin curve.

\noindent Two families of spin curves $f: \mathcal{X} \to S$ and 
$f': \mathcal{X}' \to S$ are \emph{isomorphic} if there are isomorphisms
$\sigma: \mathcal{X} \longrightarrow \mathcal{X}'$ and 
$\tau: \sigma^*(\zeta_{f'}) \longrightarrow \zeta_f$
such that $f = f' \circ \sigma$ and $\tau$ is compatible with the natural 
isomorphism between $\sigma^*(\omega_{f'})$ and $\omega_f$ (see 
\cite{Cornalba:91} p.~212). 

\noindent Let $\sgfun$ be the contravariant functor from schemes to sets, 
which to every scheme $S$ associates the set $\sgfun(S)$ of isomorphism 
classes of families of spin curves of genus $g$. 

\noindent Let $\sgbar$ be the set of isomorphism classes of spin curves of 
genus $g$ and $S_g$ be the subset consisting of classes of smooth curves. 
One can define a natural structure of analytic variety on $\sgbar$ 
(see \cite{Cornalba:89}, \S 5) and from \cite{Cornalba:89}, 
Proposition~(4.6), it follows that $\sgbar$ is a coarse moduli variety 
for $\sgfun$.

Let now $g \ge 3$. For every integer $d$, there is a  
\emph{universal Picard variety}
$$
\psi_d: \pdg \longrightarrow \mgo
$$
whose fiber $J^d(X)$ over a point $X$ of $\mgo$ parametrizes line 
bundles on $X$ of degree $d$ modulo isomorphism.

\noindent Assume $d \ge 20(g-1)$, but notice that this is not a real 
restriction because of the natural isomorphism 
$P_{d,g} \cong P_{d+n(2g-2),g}$.  
Then $\pdg$ has a natural compactification
$$
\phi_d: \pdgbar \longrightarrow \mgbar
$$ 
such that $\phi_d^{-1}(\mgo)=\pdg$.
Namely, let $\Hilb^{dx-g+1}_{d-g}$ be the Hilbert scheme parametrizing 
closed subschemes of $\P^{d-g}$ having Hilbert polynomial $d x - g + 1$ 
and set $H_d := \{ \hspace{0.1cm} 
h \in \Hilb^{dx-g+1}_{d-g}: h \textrm{\hspace{0.1cm} is \hspace{0.1cm}} 
G$-semistable and the corresponding curve is connected $\}$.
Then $\pdgbar$ was constructed in \cite{Caporaso:94} as a GIT quotient
$$
\pi_d: H_d \longrightarrow H_d / G = \pdgbar,
$$
where $G = SL(d-g+1)$.
  
\noindent Moreover, one can define (see \cite{Caporaso:94}, \S 8.1) the 
contravariant functor $\pdgfun$ from schemes to sets, which to every scheme 
$S$ associates the set $\pdgfun(S)$ of equivalence classes of polarized 
families of quasistable curves of genus $g$ 
$$
f: (\mathcal{X}, \mathcal{L}) \longrightarrow S
$$
such that $\mathcal{L}$ is a relatively very ample line bundle of degree 
$d$ whose multidegree satisfies the following Basic Inequality on each 
fiber. 

\begin{Definition}
Let $X = \bigcup_{i=1}^n X_i$ be a projective, nodal, connected curve of 
arithmetic genus $g$, where the $X_i$'s are the irreducible components of $X$. 
We say that the multidegree $(d_1, \ldots, d_n)$ 
satisfies the Basic Inequality if for every complete subcurve $Y$ of $X$ 
of arithmetic genus $g_Y$ we have
$$
m_Y \le d_Y \le m_Y + k_Y
$$
where
$$
d_Y = \sum_{X_i \subseteq Y} d_i
$$
$$
k_Y= \vert Y \cap \overline{X \setminus Y} \vert
$$
$$
m_Y= \frac{d}{g-1} \left( g_Y-1+ \frac{k_Y}{2} \right) - \frac{k_Y}{2}
$$
(see \cite{Caporaso:94} p.~611 and p.~614).
\end{Definition}

\noindent Two families over $S$, $(\mathcal{X}, \mathcal{L})$ and 
$(\mathcal{X}', \mathcal{L}')$ are \emph{equivalent} if there exists 
an $S$-isomorphism 
$$
\sigma: \mathcal{X} \longrightarrow \mathcal{X}'
$$
and a line bundle $M$ on $S$ such that
$$
\sigma^*\mathcal{L}' \cong \mathcal{L} \otimes f^*M.
$$
By \cite{Caporaso:94}, Proposition~8.1, there is a morphism of functors:
\begin{equation}\label{morphism}
\pdgfun \longrightarrow \emph{Hom} (\hspace{0.1cm} . \hspace{0.1cm}, \pdgbar)
\end{equation}
and $\pdgbar$ coarsely represents $\pdgfun$ if and only if 
\begin{equation}\label{coprime}
(d-g+1, 2g-2)=1.
\end{equation}

The relationship between $\sgbar$ and $\pdgbar$ can be expressed as follows.
\begin{Theorem}\label{naturalmorphism}
For every integer $t \ge 10$ there is a natural morphism:
$$
f_t: \sgbar \longrightarrow \overline{P}_{(2t+1)(g-1),g}.
$$
\end{Theorem}

\proof First of all, notice that in this case (\ref{coprime}) does not hold, 
so the points of $\pdgbar$ are \emph{not} in one-to-one correspondence with 
equivalence classes of very ample line bundles of degree $d$ on quasistable 
curves, satisfying the Basic Inequality (see \cite{Caporaso:94}, p.~654). 
However, we claim that the thesis can be deduced from the existence of a 
morphism of functors:
\begin{equation}\label{functor}
F_t: \sgfun \longrightarrow \overline{\mathcal{P}}_{(2t+1)(g-1),g}.
\end{equation}
Indeed, since $\sgbar$ coarsely represents $\sgfun$, any morphism of 
functors $\sgfun \to \emph{Hom} (\hspace{0.1cm} . \hspace{0.1cm}, S)$
induces a morphism of schemes $\sgbar \to S$, so the claim follows from 
(\ref{morphism}).
Now, a morphism of functors as (\ref{functor}) is the datum for any scheme 
$S$ of a set-theoretical map 
$$
F_t(S): \sgfun(S) \longrightarrow \overline{\mathcal{P}}_{(2t+1)(g-1),g}(S),
$$
satisfying obvious compatibility conditions.
Let us define
$$
F_t(S) \left( \left[ f: \mathcal{X} \to S, \zeta_f, \alpha_f \right] 
\right) := \left[ f: (\mathcal{X}, \zeta_f \otimes \omega_f^{\otimes t}) 
\to S \right].
$$
In order to prove that $F_t(S)$ is well-defined, the only non-trivial
matter is to check that the multidegree of $\zeta_f \otimes \omega_f^{\otimes 
t}$ satisfies the Basic Inequality on each fiber, so the thesis follows 
from the next Lemma.

\qed 

\begin{Lemma}\label{basic}
If $Y$ is a complete subcurve of $X$ and $d_Y$ is the degree of 
$\zeta_X \otimes \omega_X^{\otimes t} \vert_Y$, then 
$m_Y \le d_Y \le m_Y + k_Y$ in the notation of the Basic Inequality. 
Moreover, if we set  
$$
\tilde{k_Y} := \vert \tilde{Y} \cap \overline{\tilde{X} \setminus \tilde{Y}} 
\vert,
$$ 
we have $d_Y = m_Y$ if and only if $\tilde{k_Y}=0$ and all exceptional 
components in $Y$ do not intersect $\overline{X \setminus Y}$, and
$d_Y = m_Y + k_Y$ if and only if $\tilde{k_Y} = 0$ and all exceptional 
components in $\overline{X \setminus Y}$ do not intersect $Y$.
\end{Lemma}

\proof Let $Y_1, \ldots, Y_\nu$ be the irreducible components of $Y$, 
of arithmetic genus $g_1, \ldots, g_\nu$ respectively. 
We may assume that the first $\tilde{\nu}$ ones are non-ex\-ceptional and the 
last $(\nu - \tilde{\nu})$ ones are exceptional, so that 
$\tilde{Y} = Y_1 \cup \ldots \cup Y_{\tilde{\nu}}$.
Next, let $\{p_1, \ldots p_\delta \}$ be the points of intersection between 
two distinct irreducible components of $Y$. 
Again, we may assume that the first $\tilde{\delta}$ ones involve 
two non-exceptional components and the last $(\delta - \tilde{\delta})$ ones 
are between a non-exceptional and an exceptional component. We have
$$
g_Y = \sum_{i=1}^{\nu} g_i + \delta - \nu + 1 = \sum_{i=1}^{\tilde{\nu}} g_i + 
\delta - \nu + 1
$$
and since $\zeta_X^{\otimes 2} \vert_{\tilde{Y}}
\cong \omega_{\tilde{X}} \vert_{\tilde{Y}}$ we may compute 
$$
\deg \zeta_X \vert_{\tilde{Y}} = \frac{1}{2} \deg \omega_{\tilde{X}} 
\vert_{\tilde{Y}}
= \frac{1}{2} \left(\sum_{i=1}^{\tilde{\nu}} (2g_i - 2) + 2 \tilde{\delta} 
+ \tilde{k_Y} \right).
$$
Hence we deduce 
\begin{eqnarray*}
d_Y &=& \deg \left(\zeta_X \otimes \omega_X^{\otimes t} \right) \vert_Y =
\deg \zeta_X \vert_Y + t \deg \omega_X \vert_Y =\\
&=& \deg \zeta_X \vert_{\tilde{Y}} + \deg \zeta_X \vert_{\overline{Y \setminus 
\tilde{Y}}} 
+ t \deg \omega_X \vert_Y =\\
&=& \frac{1}{2} \left( \sum_{i=1}^{\tilde{\nu}} (2 g_i - 2) + 2 \tilde{\delta} 
+ \tilde{k_Y} \right) + (\nu - \tilde{\nu}) + t (2 g_Y - 2 + k_Y) =\\
&=& g_Y - 1 - (\delta - \tilde{\delta}) + 2 (\nu - \tilde{\nu}) + 
\frac{\tilde{k_Y}}{2} + 2t \left( g_Y - 1 + \frac{k_Y}{2} \right).
\end{eqnarray*}
On the other hand,
$$
m_Y = (2t+1) \left( g_Y - 1 + \frac{k_Y}{2} \right) - \frac{k_Y}{2}
= 2t \left( g_Y - 1 + \frac{k_Y}{2} \right) + g_Y -1,
$$ 
so 
$$
d_Y = m_Y - (\delta - \tilde{\delta}) + 2(\nu - \tilde{\nu}) + 
\frac{\tilde{k_Y}}{2}
$$ 
and the Basic Inequality is satisfied if and only if 
$$ 
0 \le 2(\nu -\tilde{\nu}) - (\delta - \tilde{\delta}) + \frac{\tilde{k_Y}}{2} 
\le k_Y.
$$
Now, since every exceptional component meets the other components in exactly 
two points, there are obvious inequalities
$$
\delta - \tilde{\delta} \le 2 (\nu - \tilde{\nu})
$$
and
$$
(\delta - \tilde{\delta}) + (k_Y - \tilde{k_Y}) \ge 2(\nu - \tilde{\nu}),
$$  
hence the thesis follows. 

\qed

\begin{Remark} If $t_1$ and $t_2$ are integers $\ge 10$, then
Lemma~8.1 of \cite{Caporaso:94} yields an isomorphism 
$\tau: \overline{P}_{(2t_1+1)(g-1),g} \longrightarrow 
\overline{P}_{(2t_2+1)(g-1),g}$.
We point out that by the definitions of $\tau$ (see \cite{Caporaso:94}, 
proof of Lemma~8.1) and of $f_t$ (see proof of Theorem~\ref{naturalmorphism})
there is a commutative diagram: 
$$
\sgbar \stackrel{f_{t_1}}{\longrightarrow} \overline{P}_{(2t_1+1)(g-1),g}
$$
$$
\Arrowvert 
\hspace{1.4cm} \downarrow \tau
\hspace{1.3cm}
$$
$$
\sgbar \stackrel{f_{t_2}}{\longrightarrow} \overline{P}_{(2t_2+1)(g-1),g}.
$$
\qed
\end{Remark}

\section{Spin curves in $\pdgbar$}\label{spin}
For every integer $t \ge 10$ we define
\begin{eqnarray*}
K_{(2t+1)(g-1)} := \{ h \in \Hilb^{(2t+1)(g-1)x-g+1}_{(2t+1)(g-1)-g}: 
\hspace{0.1cm} \textrm{there is a spin curve} \\ 
\hspace{0.2cm} (X, \zeta_X, \alpha_X) 
\hspace{0.2cm} \textrm{and an embedding}
\hspace{0.2cm} h_t: X \to \P^{(2t+1)(g-1)-g} \\
\hspace{0.2cm} \textrm{induced by}
\hspace{0.2cm} \zeta_X \otimes \omega_X^{\otimes t}
\hspace{0.2cm} \textrm{such that}
\hspace{0.2cm} h = h_t(X) \}.
\end{eqnarray*}
By applying \cite{Caporaso:94}, Proposition~6.1, from the first part of 
Lemma~\ref{basic} we deduce 
$$
K_{(2t+1)(g-1)} \subset H_{(2t+1)(g-1)}
$$
(the definition of $H_d$ was recalled above in section~\ref{notation}).
Moreover, the second part of the same Lemma provides a great amount of 
information on Hilbert points corresponding to spin curves.

\begin{Proposition}\label{stable}
Let $(X, \zeta_X, \alpha_X)$ be a spin curve. Then $h_t(X)$ is GIT-stable 
if and only if $\tilde{X}$ is connected.  
\end{Proposition}

\proof We are going to apply the stability criterion of \cite{Caporaso:94}, 
Lemma~6.1, which says that $h_t(X)$ is GIT-stable if and only if the only 
subcurves $Y$ of $X$ such that $d_Y = m_Y + k_Y$ are union of exceptional 
components. 

\noindent If $\tilde{X}$ is connected, then for every subcurve $Y$ of $X$ 
which is not union of exceptional components we have $\tilde{k_Y} > 0$, so 
from Lemma~\ref{basic} it follows that $d_Y < m_Y + k_Y$ and $h_t(X)$ turns 
out to be GIT-stable.

\noindent If instead $\tilde{X}$ is not connected, pick any connected 
component $Z$ of $\tilde{X}$ and take $Y$ to be the union of $Z$ with all 
exceptional components of $X$ intersecting $Z$. It follows that 
$\tilde{k_Y}=0$ and all exceptional components in $\overline{X \setminus Y}$ 
do not intersect $Y$, so Lemma~\ref{basic} yields $d_Y = m_Y + k_Y$
and $h_t(X)$ is not GIT-stable.

\qed

\begin{Proposition}\label{closed}
If $(X, \zeta_X, \alpha_X)$ is a spin curve,
then the orbit of $h_t(X)$ is closed in the semistable locus. 
\end{Proposition}

\proof Just recall the first part of \cite{Caporaso:94}, Lemma~6.1, which 
says that the orbit of $h_t(X)$ is closed in the semistable locus if and only 
if $\tilde{k_Y}=0$ for every subcurve $Y$ of $X$ such that $d_Y = m_Y$, so 
the thesis is a direct consequence of Lemma~\ref{basic}.  

\qed

The sublocus of $\overline{P}_{(2t+1)(g-1),g}$ 
obtained by projection from $K_{(2t+1)(g-1)}$ 
is indeed the GIT analogue of $\sgbar$ we are looking for.
Namely, if we set
$$
\Sigma_t := \pi_{(2t+1)(g-1)}(K_{(2t+1)(g-1)})
$$
then the following holds.

\begin{Theorem}\label{injection}
The morphism $f_t$ is an injection:
$$
f_t: \sgbar \hookrightarrow \Sigma_t.
$$
\end{Theorem}
 
\proof It is easy to check that $f_t(\sgbar)=\Sigma_t$. Indeed,
if $[(X, \zeta_X, \alpha_X)] \in \sgbar$, then any choice of a base 
for $H^0(X, \zeta_X \otimes \omega_X^{\otimes t})$ induces an 
embedding $h_t: X \to \P^{(2t+1)(g-1)-g}$ and 
$f_t([(X', \zeta_{X'}, \alpha_{X'})])
= \pi_{(2t+1)(g-1)}(h_t(X)) \in \Sigma_t$;
conversely, if $\pi_{(2t+1)(g-1)}(h) \in \Sigma_t$, then there is a 
spin curve $(X, \zeta_X, \alpha_X)$ and an embedding 
$h_t: X \to \P^{(2t+1)(g-1)-g}$ such that 
$h=h_t(X)$ and $f_t([(X', \zeta_{X'}, \alpha_{X'})])= \pi_{(2t+1)(g-1)}(h)$.

\noindent
Next we claim that $f_t$ is injective. Indeed, let $(X, \zeta_X, \alpha_X)$ 
and $(X', \zeta_{X'}, \alpha_{X'})$ be two spin curves and assume that 
$f_t([ (X, \zeta_X, \alpha_X)]) = f_t([(X', \zeta_{X'}, \alpha_{X'})])$.
Choose bases for $H^0(X, \zeta_X \otimes \omega_X^{\otimes t})$ and
$H^0(X', \zeta_{X'} \otimes \omega_{X'}^{\otimes t})$ and 
embed $X$ and $X'$ in $\P^{(2t+1)(g-1)-g}$.
If $h_t(X)$ and $h_t(X')$ are the corresponding Hilbert points, 
then $\pi_{(2t+1)(g-1)}(h(X))= \pi_{(2t+1)(g-1)}(h(X'))$ and
the Fundamental Theorem of GIT implies that 
$\overline{O_G(h_t(X))}$ and $\overline{O_G(h_t(X'))}$ intersect in the 
semistable locus. It follows from Proposition~\ref{closed} that 
$O_G(h_t(X)) \cap O_G(h_t(X')) \ne \emptyset$, so
$O_G(h_t(X))=O_G(h_t(X'))$ and there are isomorphisms
$\sigma: X \to X'$ and $\tau: \sigma^*(\zeta_{X'}) \to \zeta_X$.
Now the claim follows from Lemma~\ref{isomorphic} and the proof 
is over.

\qed

Next we are going to derive an explicit combinatorial description 
of $\Sigma_t$.
We omit the proof of the following easy Lemma, referring to the 
proof of Lemma~\ref{basic} for a similar computation.

\begin{Lemma}\label{multidegree}
Let $(X, \zeta_X, \alpha_X)$ be a spin curve. Fix a decomposition 
$$
X= \bigcup_{i=1}^n \tilde{X}_i \cup \bigcup_{j=1}^m E_j
$$
where the $\tilde{X}_i$'s are the irreducible components of $\tilde{X}$
and the $E_j$'s are the exceptional components of $X$. 
Set 
$k_i := \vert \tilde{X}_i \cap \overline{X \setminus \tilde{X}_i} \vert$ 
and 
$\tilde{k_i} := \vert \tilde{X}_i \cap \overline{\tilde{X} \setminus 
\tilde{X}_i} \vert$.
Then the multidegree of $\zeta_X \otimes \omega_X^{\otimes t}$ on $X$ 
is $\underline{d} = (d_1, \ldots, d_n, d_{n+1}, \ldots, d_{n+m})$
with
\begin{eqnarray*}
d_i &=& (2t+1) (p_a (\tilde{X}_i) - 1) + t k_i + \frac{1}{2} \tilde{k_i}  
\hspace{1cm} 1 \le i \le n \\
d_{n+j} &=& 1 \hspace{6.52cm} 1 \le j \le m.
\end{eqnarray*}
\end{Lemma}

\qed

As in \cite{Caporaso:94}, $\S$ 5.1, we set
$$
M_C^{\underline{d}} := \{ h \in H_d, h = \hilb (C,L): \underline{\deg} L
= \underline{d} \}
$$
and 
$$
V_C^{\underline{d}} := \overline{M_C^{\underline{d}}} \cap H_d.
$$
By \cite{Caporaso:94}, Corollary~5.1 (but see also on p.~627), if $[C] \in 
\mgbar$ then the $V_C^{\underline{d}}$'s are exactly the irreducible 
components of the fiber over $[C]$ of the natural morphism
$$
\psi_d: H_d \longrightarrow \mgbar.
$$

\begin{Proposition}\label{components}
Let $\underline{d} = (d_1, \ldots, d_n)$ be a multidegree and let 
$C = \bigcup_{i=1}^{n} C_i$ be a stable curve, where the $C_i$'s are 
the irreducible components of $C$. 
Set $g_i := p_a(C_i)$ and 
$k_{ij} := \vert C_i \cap C_j \vert$ if $i \ne j$, $k_{ij} := 0$
if $i = j$. 
\newline \noindent
Then there exists a spin curve $(X, \zeta_X, \alpha_X)$ such that 
$h_t(X) \in  V_C^{\underline{d}}$ if and only if for every $1 \le i,j \le n$ 
there are integers $s_{ij}$ and $\sigma_{ij}$ with
\begin{eqnarray*}
0 \le &s_{ij}& \le k_{ij} \hspace{1cm}
s_{ij} = s_{ji} \hspace{1cm}
\sum_{j=1}^n (k_{ij} - s_{ij}) \equiv 0  \hspace{0.1cm} \mod 
\hspace{0.1cm} 2 \\
0 \le &\sigma_{ij}& \le s_{ij} \hspace{1cm}
\sigma_{ij} + \sigma_{ji} = s_{ij} 
\end{eqnarray*}
such that
$$
d_i = (2t+1) (g_i - 1) + t \sum_{j=1}^n k_{ij} 
+ \frac{1}{2} \sum_{j=1}^n (k_{ij} - s_{ij})  
+ \sum_{j=1}^n \sigma_{ij}.
$$
\end{Proposition}    

\proof To get a quasistable curve $X$ starting from 
$C = \bigcup_{i=1}^{n} C_i$, for every $1 \le i,j \le n$
choose $r_i$ nodes of $C_i$ and 
$s_{ij}$ contact points between $C_i$ and $C_j$
and blow them up. 
In the notation of Lemma~\ref{multidegree}, notice that 
$p_a(\tilde{X}_i) = g_i - r_i$, $k_i = \sum_{j=1}^n k_{ij} + 2 r_i$
and $\tilde{k_i} = \sum_{j=1}^n (k_{ij} - s_{ij})$.
As pointed out in \cite{Cornalba:89}, $\S$ 3, in order for a spin curve 
having $X$ as underlying curve to exist, a necessary and sufficient 
condition is that $\tilde{k_i} \equiv 0  \hspace{0.1cm} \mod 
\hspace{0.1cm} 2$ for every $1 \le i \le n$. Moreover, by 
\cite{Caporaso:94}, Proposition~5.1, $h_t(X) \in V_C^{\underline{d}}$
if and only if there is a partition $X = \bigcup_{i=1}^n X_i$ such that
$X_i$ is a complete connected subcurve of $X$ whose stable model is $C_i$ 
and $d_i = \deg_{X_i} (\zeta_X \otimes \omega_X^{\otimes t})$. 
So the thesis follows from Lemma~\ref{multidegree}.

\qed

Let $V_C := \bigcup_{\underline{d}} V^{\underline{d}}_C$ and let 
$\pdC := \phi_d^{-1}(C)$.  
By \cite{Caporaso:94}, proof of Corollary~5.1, we have
$$
\pdC = V_C / G
$$
and for every irreducible component $I$ of $\pdC$ there is a unique 
multidegree $\underline{d}$ such that $V^{\underline{d}}_C$ dominates 
$I$ via the quotient map 
$$
V_C \longrightarrow \pdC.
$$ 

\begin{Theorem}\label{image}
Let $C = \bigcup_{i=1}^{n} C_i$ be a stable curve, where the $C_i$'s are 
the irreducible components of $C$. 
Set $g_i := p_a(C_i)$ and 
$k_{ij} := \vert C_i \cap C_j \vert$ if $i \ne j$, $k_{ij} := 0$
if $i = j$. 
Let $I$ be an irreducible component of $\pdC$ and let $\underline{d} =
(d_1, \ldots, d_n)$ be the multidegree such that $I$ is dominated by 
$V^{\underline{d}}_C$.
\newline \noindent
Then there exists a spin curve $(X, \zeta_X, \alpha_X)$ such that 
$f_t ( \left[ (X, \zeta_X, \alpha_X) \right] ) \in I$ if and only if 
for every $1 \le i,j \le n$ 
there are integers $s_{ij}$ and $\sigma_{ij}$ with
\begin{eqnarray*}
0 \le &s_{ij}& \le k_{ij} \hspace{1cm}
s_{ij} = s_{ji} \hspace{1cm}
\sum_{j=1}^n (k_{ij} - s_{ij}) \equiv 0  \hspace{0.1cm} \mod 
\hspace{0.1cm} 2 \\
0 \le &\sigma_{ij}& \le s_{ij} \hspace{1cm}
\sigma_{ij} + \sigma_{ji} = s_{ij} 
\end{eqnarray*}
such that
$$
d_i = (2t+1) (g_i - 1) + t \sum_{j=1}^n k_{ij} 
+ \frac{1}{2} \sum_{j=1}^n (k_{ij} - s_{ij})  
+ \sum_{j=1}^n \sigma_{ij}.
$$
\end{Theorem}

\proof By the Fundamental Theorem of GIT, 
$$
f_t([(X, \zeta_X, \alpha_X)]) = \pi_{(2t+1)(g-1)}(h_t(X)) \in I
$$
if and only if there is $h \in V^{\underline{d}}_C$ such that 
$\overline{O_G(h_t(X))}$ and $\overline{O_G(h)}$ intersect in the 
semistable locus. 

\noindent Since $O_G(h_t(X))$ is closed in the semistable locus by 
Proposition~\ref{closed}, we have
$$
O_G(h_t(X)) \cap \overline{O_G(h)} \ne \emptyset
$$
and since $\overline{O_G(h)}$ is a union of orbits we may rephrase the 
above condition as
$$
h_t(X) \in \overline{O_G(h)}.
$$
On the other hand, we have
$V^{\underline{d}}_C = \bigcup_{h \in V^{\underline{d}}_C} O_G(h)$
since $V^{\underline{d}}_C$ is $G$-invariant and 
$V^{\underline{d}}_C = \bigcup_{h \in V^{\underline{d}}_C} 
\overline{O_G(h)}$
since $V^{\underline{d}}_C$ is closed. 

\noindent Summing up, we see that 
$f_t ( \left[ (X, \zeta_X, \alpha_X) \right] ) \in I$ 
if and only if $h_t(X) \in V^{\underline{d}}_C$.
Now the thesis follows from Proposition~\ref{components}.

\qed 

\noindent
\textbf{Example.} Let $C$ be a \emph{split curve} of genus $g$, i.e. the 
union of two nonsingular rational curves meeting transversally at $g+1$ 
points. Such curves are particularly interesting, for a number of reasons 
(see \cite{Caporaso:01} and \cite{CapSer:02}). 
According to Theorem~\ref{image}, $\Sigma_t$ meets an irreducible component 
$I$ of $\pdC$ if and only if $I$ corresponds to a bidegree $(d_1, d_2)$
with 
\begin{eqnarray*}
d_1 &=& \left(t + \frac{1}{2} \right) (g+1) - (2t+1) - \frac{1}{2} s + 
\sigma \\
d_2 &=& \left(t + \frac{1}{2} \right) (g+1) - (2t+1) + \frac{1}{2} s - 
\sigma
\end{eqnarray*}
where $s$ and $\sigma$ are nonnegative integers satisfying
$$
\sigma \le s \le g+1 \hspace{2cm} s \equiv g+1  \hspace{0.1cm} \mod 
\hspace{0.1cm} 2.
$$
Moreover, from the proof of Proposition~\ref{components} it follows that 
$s$ is the number of exceptional components of a quasi-stable curve $X$ 
underlying a spin curve $(X, \zeta_X, \alpha_X)$ such that 
$f_t([(X, \zeta_X, \alpha_X)]) \in \Sigma_t \cap I$.
\qed

\section{Divisors on $\pdgbar$}\label{divisors}
In order to understand the boundary of $\pdgbar$, we recall the  
decomposition of the boundary of $\mgbar$ into its irreducible components:
$$
\partial \mgbar = \Delta_0 \cup \Delta_1 \cup \ldots \cup \Delta_{\lfloor 
g / 2 \rfloor}
$$
and we define
$$
D_i := \phi_d^{-1}(\Delta_i)
$$
for $i=0, \ldots \lfloor g / 2 \rfloor$.
Notice that, since $\phi_d$ is surjective, each $D_i$ turns out to be 
a divisor on $\pdgbar$.
Moreover, if $X \in \mgbar$, we set as usual $\pdX := \phi_d^{-1}(X)$.
 
\begin{Lemma}\label{irreducible} For every $i$, if $X$ is a general 
element in $\Delta_i$ then $\pdX$ is irreducible. 
\end{Lemma}

\proof If $i \ge 1$, a general element $X$ of $\Delta_i$ is the union of 
two smooth curves $X_1$ and $X_2$ meeting at one node, so $X$ is of 
compact type and $\pdX$ is irreducible (see \cite{Caporaso:94}, footnote 
on p.~594). If instead $i=0$, a general element 
$X \in \Delta_0$ is an irreducible curve and also $\pdX$ turns out to be 
irreducible (see \cite{Caporaso:94} 7.1).

\qed

As a consequence, we obtain a complete description of the boundary of 
$\pdgbar$.
\begin{Proposition}\label{unique} For every $i$, $D_i$ is irreducible.
\end{Proposition}

\proof By \cite{Caporaso:94}, Corollary~5.1 (2), for every $X \in \mgbar$ 
all irreducible components of $\phi_d^{-1}(X)$ have dimension $g$. 
Let $I$ be an irreducible component of $D_i$; by applying 
the theorem on the dimensions of the fibers (see for instance 
\cite{Hartshorne:77}, II, Ex.~3.22 (b) p.~95) to the map
$$
\phi_{d \vert I}: I \longrightarrow \phi_d(I) \subseteq \Delta_i
$$
we obtain
$$
\dim I - \dim \phi_d(I) \le \dim (\phi_d^{-1}(X) \cap I) \le g.
$$
Hence
\begin{eqnarray*}
3g-4 &=& \dim \Delta_i \ge \dim \phi_d(I) \ge 4g-4 - \dim(\phi_d^{-1}(X) 
\cap I) \\
&\ge& 3g-4
\end{eqnarray*}
and all the above inequalities turn out to be equalities. In particular, 
we have
$$
\dim(\phi_d^{-1}(X) \cap I) = g
$$
and $\phi_d^{-1}(X) \cap I$ is a closed subscheme of $\phi_d^{-1}(X)$ of 
maximal dimension. By Lemma~\ref{irreducible} there is a dense open 
subset $U_i \subset \Delta_i$ such that $\phi_d^{-1}(X)$ is irreducible 
for every $X \in U_i$. It follows that $\phi_d^{-1}(X) \subset I$ for 
every $X \in U_i$ and  $\phi_d^{-1}(U_i)$ is a dense open subset of $I$. 
Hence $I = \overline{\phi_d^{-1}(U_i)}$ is uniquely determined.

\qed

We recall that the class of any line bundle $\mathcal{L}$ on $\pdg$
restricted to a fiber $J^d(X)$ is a multiple $m \theta$ of the class $\theta$ 
of the $\Theta$ divisor (see \cite{Kouvidakis:91} p.~840); it seems therefore 
natural to define the \emph{class} of $\mathcal{L}$ to be the integer $m$.

The Picard group of $\pdg$ is completely described by the following 
result, due to A.~Kouvidakis (see \cite{Kouvidakis:91}, Theorem~4 p.~849). 
\begin{Theorem}~\label{Kouvidakis} If $g \ge 3$ then 
the Picard group of the universal Picard variety $\psi_d: \pdg \to \mgo$ 
is freely generated over $\Q$ by the line bundles  $\mathcal{L}_{d,g}$ 
and $\psi_d^{*}(\lambda)$, where $\mathcal{L}_{d,g}$ is any line bundle 
on $\pdg$ with class 
$$
k_{d,g}= \frac{2g-2}{\gcd (2g-2, g+d-1)}
$$
and $\lambda$ is the Hodge bundle on $\mgo$.
\end{Theorem}
 
The analogous result for $\pdgbar$ is the following:

\begin{Theorem}\label{A_4g-4} Assume $g \ge 3$ and $d \ge 20(g-1)$. 
Then the divisor class group of the universal Picard variety $\phi_d: 
\pdgbar \to \mgbar$ is freely generated over $\Q$ by the classes 
$\mathcal{L}_{d,g}$, $\phi_d^{*}(\lambda)$ and 
$\mathcal{D}_i$ ($i=0, \ldots \lfloor g / 2 \rfloor$), where  
$\mathcal{D}_i$ denotes the linear equivalence class of $D_i$.
\end{Theorem}

\proof Since $\pdg$ is smooth and irreducible, there is a natural 
identification $\Pic(\pdg)=A_{4g-4}(\pdg)$, so using the exact sequence
$$
A_{4g-4}(\pdgbar \setminus \pdg) \to A_{4g-4}(\pdgbar) \to A_{4g-4}(\pdg)
\to 0
$$
we may deduce from Theorem~\ref{Kouvidakis} that $A_{4g-4}(\pdgbar)$ is 
generated by $\mathcal{L}_{d,g}$, $\phi_d^{*}(\lambda)$ and 
$A_{4g-4}(\pdgbar \setminus \pdg)$. Now, we have
$$
\pdgbar \setminus \pdg = D_0 \cup \ldots \cup D_{\lfloor g / 2 \rfloor} 
\cup \phi_d^{-1}(\mg \setminus \mgo).
$$
Moreover, by applying the theorem on the dimensions of the fibers to the map
$$
\phi_{d \vert \phi_d^{-1}(\mg \setminus \mgo)}: \phi_d^{-1}(\mg \setminus \mgo)
\longrightarrow \mg \setminus \mgo
$$
we obtain that $\codim(\phi_d^{-1}(\mg \setminus \mgo), \pdgbar) 
\ge \codim(\mg \setminus \mgo, \mgbar) = g-2$. 
Hence, if $g \ge 4$, Proposition~\ref{unique} implies that $A_{4g-4}(\pdgbar 
\setminus \pdg)$ is generated by the $\mathcal{D}_i$'s. 
If, instead, $g=3$, we recall that the hyperelliptic locus $H$ is the unique 
divisor in $\mathcal{M}_3$ contained in $\mathcal{M}_3 \setminus 
\mathcal{M}_3^0$ (see \cite{HarMor:98}, Ex.~2.27,~3)). Since 
$\left[ H \right] = 18 \lambda$ in $\Pic(\mathcal{M}_3 \otimes \Q)$
(see \cite{HarMor:98} p.~164), we have $\left[ \phi_d^{-1}(H) \right] = 
18 \phi_d^*(\lambda)$ and the result on generation is completely proved. 
As for relations, let 
\begin{equation}\label{rel}
a \mathcal{L}_{d,g} + b \phi_d^{*}(\lambda) + \sum c_i \mathcal{D}_i = 0
\end{equation}
be a relation in $A_{4g-4}(\pdgbar)$. If $J^d(X)$ is the fiber over a  
curve $X$ in $\mgo$, then restricting (\ref{rel}) to $J^d(X)$ yields 
$ a k_{d,g} \theta = 0$ and we get $a = 0$. So we may rephrase (\ref{rel})
as $\phi_d^*(b \lambda + \sum c_i \delta_i) = 0$. 
The thesis is now a direct consequence of the following Lemma.

\qed 

\begin{Lemma}\label{injective}
Assume $g \ge 3$ and $d \ge 20(g-1)$. Then 
$$
\phi_d^*: A_{3g-4}(\mgbar) \otimes \Q \longrightarrow A_{4g-4}(\pdgbar) 
\otimes \Q
$$
is injective.
\end{Lemma}

\proof It is well-known that $A_{3g-4}(\mgbar) \otimes \Q$ is generated 
by the Hodge class $\lambda$ and the boundary classes $\delta_i$ (see 
for instance \cite{ArbCor:87}, Proposition~2). A natural way to check 
that these classes are independent is to construct families of 
stable curves $h: X \to S$, with $S$ a smooth complete curve, such that 
the vectors $(\deg_h(\lambda), \deg_h(\delta_0), \ldots)$ are linearly 
independent. Such a construction is carried out in detail in 
\cite{ArbCor:87} and it is applied in \cite{Cornalba:89} to show the 
injectivity of $\chi^{*}: \Pic(\mgbar) \to \Pic(\sgbar)$. 
We are going to mimic the same idea in our case. Namely, in order to prove 
that $\phi_d^{*}(\lambda)$ and the $\phi_d^{*}(\delta_i)$'s are independent  
in $A_{4g-4}(\pdgbar) \otimes \Q$, we will lift to $\pdgbar$ the families 
$h: X \to S$ constructed in \cite{ArbCor:87}. The key observation is that 
each of them is equipped with many sections passing through the smooth 
locus of the general curve $C$ of the family. Indeed, for every irreducible 
component $C_i$ of $C$ we easily find a section $\sigma_i$ of $h$ which cuts 
on $C$ a smooth point $P_i \in C_i$. Next, we decompose the integer $d$ as a 
sum of $d_i$'s in such a way that the multidegree determined by the $d_i$'s 
satisfies the Basic Inequality. Finally, we endow $C$ with the line bundle 
$\mathcal{L} := \otimes_i \mathcal{O}_{C_i}(d_i P_i)$. Since $\pdgbar$ is 
proper over $\mgbar$, this construction uniquely determines a lifting of 
$h: X \to S$, so the proof is over.

\qed

\begin{Corollary}\label{Pic} Assume $g \ge 3$, $d \ge 20(g-1)$ and 
$(d-g+1, 2g-2)=1$. Then $\Pic(\pdgbar) \otimes \Q$ is freely generated 
by $\mathcal{L}_{d,g}$, $\phi_d^{*}(\lambda)$ and $\mathcal{D}_i$ 
($i=0, \ldots \lfloor g / 2 \rfloor$).
\end{Corollary}

\proof By \cite{Caporaso:94}, Lemma~2.2~(1), $\pdgbar$ is the quotient of a 
nonsingular scheme, so in particular it is normal and there is an injection:
$$
\Pic(\pdgbar) \hookrightarrow A_{4g-4}(\pdgbar).
$$
Moreover, since $(d-g+1, 2g-2)=1$, the singularities of $\pdgbar$ are all of 
finite quotient type (see \cite{Caporaso:94}, Proposition on p.~594). 
It follows that every Weil divisor is $\Q$-Cartier, so we get a surjective 
morphism:
$$
\Pic(\pdgbar) \otimes \Q \twoheadrightarrow A_{4g-4}(\pdgbar) \otimes \Q.
$$
Hence Theorem~\ref{A_4g-4} yields the thesis.

\qed

\vspace{0.5cm}
\noindent
Claudio Fontanari \newline
Universit\`a degli Studi di Trento \newline
Dipartimento di Matematica \newline
Via Sommarive 14 \newline
38050 Povo (Trento) \newline
Italy \newline
e-mail: fontanar@science.unitn.it

\end{document}